\documentclass[a4paper]{amsart}
\usepackage{mathrsfs,amssymb}
\usepackage{multicol}\multicolsep=0pt
\usepackage{pstricks,pst-node}
\usepackage{cite}

\def\crulefill{\leavevmode\leaders\hrule height 1pt\hfill\kern 0pt}
\long\def\QUERY#1{%
\leavevmode\newline%
\noindent$\star\star\star$\thinspace\textsf{Comment/Query}\crulefill\newline%
   \space #1\newline\hbox to 120mm{\crulefill}$\star\star\star$\newline
}

\numberwithin{equation}{section} \theoremstyle{definition}
\newtheorem{Defn}[equation]{Definition}
\theoremstyle{plain}
\newtheorem{Prop}[equation]{Proposition}
\newtheorem{Theorem}[equation]{Theorem}
\newtheorem{Lemma}[equation]{Lemma}
\newtheorem{Cor}[equation]{Corollary}
\newtheorem{Remark}[equation]{Remark}

\newtheorem{THEOREM}{Theorem}

\makeatletter
\def\enumerate{\begingroup\ifnum\@enumdepth>3\@toodeep\else
      \advance\@enumdepth\@ne
      \edef\@enumctr{enum\romannumeral\the\@enumdepth}%
      \topsep\z@\parskip\z@
      \list{\csname label\@enumctr\endcsname}
        {\@nmbrlisttrue\let\@listctr\@enumctr
         \parsep\z@\itemsep\z@\topsep\z@
         \setcounter{\@enumctr}{0}
         \def\makelabel##1{\hss\llap{\rm ##1}}
       }\fi}

\makeatother

\let\bar=\overline
\let\epsilon=\varepsilon
\def\({\big(}
\def\){\big)}

\def\G{\mathcal G}

\newcommand{\cba}[1]{\mathscr B_{m,#1}(\boldsymbol{\delta})}
\newcommand{\ccba}[1]{\mathscr B_{m,#1}(\boldsymbol{0})}

\def\H{\mathscr H}

\DeclareMathOperator{\Char}{char}

\DeclareMathOperator{\Rad}{Rad}

\def\la{\lambda}
\def\La{\Lambda}

\def\s{\mathfrak s}

\def\t{\mathfrak t}

\def\ind{\text{Ind}}

\newcommand{\iocba}[2]{I_{m,#1}^{#2 '}(\boldsymbol{0})}

{\catcode`\|=\active
  \gdef\set#1{\mathinner{\lbrace\,{\mathcode`\|"8000%
                                   \let|\midvert #1}\,\rbrace}}
}
\def\midvert{\egroup\mid\bgroup}

\def\Set[#1]#2|#3|{\Big\{\ #2\ \Big| \
           \vcenter{\hsize #1mm\centering #3}\Big\}}

\title{On the semisimplicity of the cyclotomic Brauer algebras, II}

\author{Hebing Rui}
\address{Department of Mathematics, East China Normal University, %
         200062 Shanghai, P.R. China.}
\email{hbrui@math.ecnu.edu.cn}

\author{Jie Xu}
\address{Department of Mathematics, East China Normal University, %
         200062 Shanghai, P.R. China.}
\email{52060601009@student.ecnu.edu.cn}
\thanks{The first author was supported by NSFC No. 10331030 and NCET}

\begin{document}
\baselineskip15pt
\begin{abstract} In this paper, we give a necessary and sufficient
condition for the semisimplicity of  cyclotomic Brauer algebras
$\cba{n}$ of types $G(m, 1, n)$ with $m\ge 2$.  This generalizes
\cite[1.2--1.3]{R} and \cite[2.5]{RS} on Brauer algebras.
\end{abstract}
 \sloppy \maketitle

\centerline{\textit{Dedicated to Professor Gordon James on the
occasion of  his 60th birthday}}

\section{Introduction}

The cyclotomic Brauer algebras $\cba{n}$ have been  introduced by
H$\ddot{\text{a}}$ring-Oldenburg in \cite{O} as  classical limits of
cyclotomic Birman-Murakami-Wenzl algebras. When $m=1$, they are
Brauer algebras $\mathscr B_n(\delta)$\cite{B}.

The main purpose of this paper is to give a necessary and sufficient
condition for the semisimplicity of  $\cba{n}$ under the assumption
$m\ge 2$. For $m=1$, such a criterion has been given in
\cite[1.2-1.3]{R} and \cite[2.5]{RS}.

Unless otherwise stated, we assume that $F$ is  a splitting field of
$x^m-1$, which contains $\delta_i, 1\le i\le m$. By assumption,
there are $u_i\in F$ such that $x^m-1=\prod_{i=1}^m (x-u_i)$. Define
\begin{equation}\label{e}
e=\begin{cases}  +\infty,  & \text{if $\Char F=0$,}\\
 \Char F, & \text{if $\Char F>0$.}\\
 \end{cases}
\end{equation}

Following \cite{RS}, we define $\mathbb Z_{m, n}=\{ma\mid a\in
\tilde {\mathbb Z}_{m, n}\}$, where $\tilde {\mathbb Z}_{m, n}$ is
given as follows:

\begin{itemize} \item[(1)] $\tilde
{\mathbb Z}_{2, n}=\tilde{\mathbb Z}_{1, n}=\{k\in \mathbb Z\mid
3-n\le k\le n-3\}\cup \{2k-3\mid 3\le k\le n, k\in \mathbb Z\}$.

\item[(2)] $\tilde{\mathbb Z}_{m,n}=\tilde{\mathbb Z}_{1, n}\cup \{2-n, n-2\}$ if  $m\geq3$ and $n\ge 2$.
\end{itemize}

Suppose that  $x_1, x_2, \cdots, x_m$ are indeterminates over $F$.
If $F$ contains $\xi$, a primitive $m$-th root of unity, then we
define
\begin{equation}\label{bardelta}
\bar x_i=\sum_{j=1}^m x_j \xi^{ji},\quad   0\le i\le
m-1.\end{equation}
 Note that $F$ contains $\xi$ if   $e\nmid m$ \cite[8.2]{H}. The
following is the main result of  this paper.

\begin{THEOREM}\label{RANK} Fix two positive integers $m, n$ with $m>1$.
Let $\cba{n}$ be a cyclotomic Brauer algebra over $F$.

\begin{itemize}
\item[(a)] Suppose $n\ge 2$.
If $\delta_i\neq 0$ for some $i$, $0\le i\le m-1$, then
   $\cba{n}$ is (split) semisimple if and only if
\begin{itemize}\item[(1)] $e\nmid m\cdot n!$,
\item[(2)] $\epsilon_{i, 0} m-\bar\delta_i\not\in \mathbb Z_{m, n}$,
$0\le i\le m-1$, where $\epsilon_{i, 0}$ is the Kronecker
function.
\end{itemize}

\item[(b)] Suppose $n\ge 2$. If $\delta_i=0$, $0\le i\le m-1$, then $\ccba{n}$ is
not (split) semisimple.

\item[(c)]   $\cba{1}$ is (split) semisimple if and only if  $e\nmid m$.
\end{itemize}
\end{THEOREM}

In what follows, we write $\delta_j=\delta_i$ if $i, j\in \mathbb Z$
and $i\equiv j \mod m$.

Let $\H_{i, k}$ be the hyperplane in $F^m$, which is determined by
the linear function $\epsilon_{i, 0} m-\bar x_i=k$, $0\le i\le m-1$
and $k\in \mathbb Z_{m, n}$.  Condition (2) in Theorem~\ref{RANK}(a)
is equivalent to the fact that $(\delta_0, \delta_1,
\cdots,\delta_{m-1})\not\in  \cup_{0\le i\le m-1, k\in \mathbb Z_{m,
n}} \H_{i, k}$. When $m=1$,  $\H_{i, k}$ collapses to a point in
$1$-dimensional $F$-space. This result has been proved in
\cite[1.2-1.3]{R} and \cite[2.5]{RS}. We remark that certain
sufficient conditions for  semisimplicity of complex Brauer algebras
have been given in \cite{Br, DWH, W}.

Our proof depends on Graham-Lehrer's theory on cellular algebras
\cite{GL} and Doran-Wales-Hanlon's work \cite[3.3-3.4]{DWH} on
Brauer algebras. Let's explain the idea as follows.

In \cite{GL}, Graham and Lehrer have introduced the notion of
cellular algebra which is defined over a poset $\La$. Such an
algebra has a nice basis, called a cellular basis. For each
$\lambda\in \Lambda$, one can define $\Delta(\lambda)$, called a
cell module. Graham and Lehrer have shown that there is a symmetric,
associative bilinear form $\phi_{\lambda}$ defined on
$\Delta(\lambda)$. It has been proved in \cite[3.8]{GL} that a
cellular algebra is (split) semisimple if and only if
$\phi_{\lambda}$ is non-degenerate for any $\lambda\in \Lambda$. It
is well known that a cellular algebra is split semisimple if and
only if it is semisimple. Therefore, one can determine  whether a
cellular algebra is semisimple by deciding  if all $\phi_{\lambda}$
are non-degenerate.

In \cite{GL}, Graham-Lehrer  have proved that a Brauer algebra
$\mathscr B_n(\delta)$ over a commutative ring is a cellular algebra
over the poset $\La$ which consists of all pairs $(f, \lambda)$,
with $0\le f\le \lfloor n/2\rfloor$ and $\lambda$ being a partition
of $n-2f$. Here $\lfloor n/2\rfloor$ is the maximal integer which is
less than $n/2$. Therefore, one can study the semisimplicity of
$\mathscr B_n(\delta)$  by deciding  whether $\phi_{f, \lambda}$ is
non-degenerate or not for any $(f,\lambda)\in \La$. Unfortunately,
it is difficult to determine whether $\phi_{f, \lambda}$ is
degenerate or not for a fixed $(f, \lambda)$.

In \cite{R}, the first author has proved that  the semisimplicity of
$\mathscr B_n(\delta)$ is completely determined by  $\phi_{f,
\lambda}$ for all partitions $\lambda$ of $n-2f$ with $f=0, 1$.
Using \cite[3.3-3.4]{DWH}, he has  decided whether such $\phi_{f,
\lambda}$'s are degenerate or not in \cite{R}. This   gives a
complete solution of the problem of   semisimplicity of $\mathscr
B_n(\delta)$   over an arbitrary  field. This method will be   used
to study the semisimplicity of $\cba{n}$ in the current paper.

The contents of this paper are organized as follows. In section~2,
we state some results on  cyclotomic Brauer algebras, and  complex
reflection group $W_{m, 2}$.  In section~3, we describe explicitly
the zero divisors of the discriminants for certain cell modules.
Theorem~A will be proved in section~4.

\textbf {Acknowledgement:} We thank the referee for his/her helpful
comments.

\section{Cyclotomic Brauer algebras}

Let $R$ be a commutative ring which contains the identity $1_R$ and
$\delta_i, 1\le i\le m$ . The cyclotomic Brauer algebra $\cba{n}$
with parameters $\delta_i, 1\le i\le m$, is the associative
$R$-algebra which is free as $R$-module with basis which consists of
all labelled Brauer diagrams \cite{O}. $\cba{n}$ can also be defined
as the $R$-algebra generated by $\set{s_i,e_i, t_j|1\le i<n \text{
and }1\le j\le n } $ subject to the relations
\begin{multicols}{2}
\begin{enumerate}
\item $s_i^2=1$, for $1\le i<n$.
\item $s_is_j=s_js_i$ if $|i-j|>1$.
\item $s_is_{i+1}s_i=s_{i+1}s_is_{i+1}$,\newline for $1\le i<n-1$.
\item $s_it_j=t_js_i$ if $j\ne i,i+1$.
\item $e_i^2=\delta_0e_i$, for $1\le i<n$.
\item $s_ie_j=e_js_i$, if $|i-j|>1$.
\item $e_ie_j=e_je_i$, if $|i-j|>1$.
\item $e_it_j=t_je_i$, if $j\ne i,i+1$.
\item $t_it_j=t_jt_i$, for $1\le i,j\le n$.
\item $s_it_i=t_{i+1}s_i$
 for $1\le i<n$.
\item $e_is_i=e_i=s_ie_i$, for $1\le i\le n-1$.
\item $s_ie_{i+1}e_i=s_{i+1}e_i$, \newline for $1\le i\le n-2$.
\item $e_{i+1}e_is_{i+1}=e_{i+1}s_i$, \newline for $1\le i\le n-2$.
\item $e_{i}e_je_{i}=e_{i}$ if $|i-j|=1$.
\item $e_it_it_{i+1}=e_i=t_it_{i+1}e_i$, \newline for $1\le i<n$.
\item $e_it_i^ae_i=\delta_a e_i$, for $1\le a\le m-1$ \newline and
$1\le i\le n-1$.
\item $t_i^m=1$, for $1\le i\le n$.
\end{enumerate}
\end{multicols}

One can prove that the two definitions of $\cba{n}$  are equivalent
by the arguments similar to those for Brauer algebras  in \cite{MW}.


The following result can be proved easily by checking the  defining
relations of $\cba{n}$.

\begin{Lemma}\label{quasi} Let $\cba{n}$ be a cyclotomic Brauer
algebra over  $R$.
There is an $R$-linear anti-involution
 \label{invo} $\ast: \cba{n}\rightarrow \cba{n}$
 such that $h^*=h$ for all $h\in \{e_i, s_i, t_j\mid 1\le i<n, 1\le j\le n\}$.



\end{Lemma}

Recall that $F$ is a splitting field of $x^m-1$. In the remaining
part of this  section, we assume  $e\nmid m\cdot n!$. By
\cite[8.2]{H}, $F$ contains $\xi$, a primitive $m$-th root of unity.

We will decompose an $FW_{m,2}$-module in Proposition~\ref{eta},
where $W_{m, n}$ is the complex reflection group of type $G(m, 1,
n)$. Note that  $W_{m, n}$ is generated by $s_i, t_1 $ satisfying
the relations
\begin{itemize}
\item $s_i^2=t_1^m=1$ for $1\le i\le n-1$.
\item $s_is_js_i=s_js_is_j$, if $|i-j|=1$.
\item $s_is_j=s_js_i$,  if $|i-j|>1$.
\item $s_it_1=t_1s_i$ if $1<i\le n-1$.
\item $s_1t_1s_1t_1=t_1s_1t_1s_1$.\end{itemize}

The order of $W_{m, n}$ is $m^n \cdot n!$. By Maschke's theorem,
the group algebra $F W_{m, n}$ is (split) semisimple.

Let $\Lambda_{m}^+(n)$ be the set of $m$-partitions of $n$. When
$m=1$, we use $\Lambda^+(n)$ instead of  $\Lambda_{1}^+(n)$. For
any $\lambda\in \Lambda_m^+(n)$, let $S^\lambda$ be the classical
Specht module with respect to $\lambda$ (see \cite[2.1]{DR1}).

For any $\lambda\in \La^+(n)$, let $\mu=(\mu_1, \mu_2, \cdots)$
with $\mu_i=\# \{j\mid \lambda_j\ge i\}$. Then $\mu$, which will
be denoted by $\lambda'$, is called the dual partition of
$\lambda$. If $\lambda=(\lambda^{(1)},  \lambda^{(2)}, \cdots,
\lambda^{(m)})\in \Lambda_m^+(n)$, we write
$\lambda'=(\lambda^{(m)'},  \lambda^{(m-1)'}, \cdots,
\lambda^{(1)'})$ and call $\lambda'$ the dual partition of
$\lambda$.

\begin{Remark}\textsf{ All  modules considered in this
paper are left modules}. I.e. $S^{\lambda}=FW_{m, n} y_{\lambda'}
w_{\lambda'} x_{\lambda} $ if we keep the notation in \cite{DR1}. In
\cite{RY}, we have assumed  $u_i=\xi^{i}$, $1\le i\le m$. In this
paper, we keep this assumption  in order to use results in \cite{RY}
directly.
\end{Remark}
Since $FW_{m, n}$ is the Ariki-Koike algebra with  $q=1$ and
$x_1^m-1=\prod_{i=1}^m (x_1-u_i)$, the following result is a special
case of the result in   \cite{DR1}.

\begin{Lemma} The set $\set { S^\lambda\mid \lambda\in
\Lambda_m^+(n)}$ is a complete set of pairwise non-isomorphic
irreducible $FW_{m, n}$-modules.\end{Lemma}

\begin{Defn} Let $m$ be a positive integer. If $m$ is even,  we define
$\wp_m(2)=\{\eta_i\mid \frac m2\le i\le m\}$, where
  $$\eta_i=\begin{cases} (0, \cdots, 0, 2), & \text{if $i=m$,}\\
(0, \cdots, 0, \underset {\frac{m}{2}-th} 2, 0\cdots, 0),
&\text{if $i=\frac{m}{2}$,}\\
(0, \cdots, 0, \underset {(m-i)-th} 1, 0, \cdots, 0, \underset
{i-th} 1, 0, \cdots, 0), &\text{if  $\frac m2< i\le
m-1$}.\\
\end{cases}
$$

If $m$ is odd, we define  $\wp_m(2)=\{\eta_i\mid \frac{m+1}2\le i\le
m\}$, where $$
\eta_i=\begin{cases} ( 0, \cdots, 0, 2), & \text{if $i=m$,}\\
(0, \cdots, 0, \underset {(m-i)-th} 1, 0, \cdots, 0, \underset
{i-th} 1, 0, \cdots, 0), & \text{if $\frac{m+1}2  \le
i\le m-1$.}\\
\end{cases}
$$
\end{Defn}

\begin{Prop}\label{eta}  Let $\mathbb Z_m\wr B_1$ be the
subgroup of $W_{m, 2}$ generated by $s_1$, $t_1t_2$.  As $FW_{m,
2}$-modules, $\ind_{\mathbb Z_m\wr B_1}^{W_{m, 2}} 1\cong
\bigoplus_{\boldsymbol \eta\in \wp_m(2)} S^{\boldsymbol \eta}$.
\end{Prop}

\begin{proof} Since  $\{1, t_1, \cdots, t_1^{m-1}\}$ is a
complete set of left coset representatives of $\mathbb Z_m\wr B_1$
in $ W_{m, 2}$, $ \{t_1^k\sum_{l=0}^{m-1} (t_1t_2)^l (1+s_1) \mid
0\le k\le m-1\}$ is an $F$-basis of $\ind_{\mathbb Z_m\wr
B_1}^{W_{m, 2}} 1$.  By assumption,  $F$ contains a primitive $m$-th
root of unity, say $\xi$. Since we are assuming that $u_i=\xi^{i}$,
$1\le i\le m$,
  $\ind_{\mathbb Z_m\wr B_1}^{\mathbb Z_m\wr \mathfrak S_2} 1$ has a
basis $\set{w_i\mid 1\le i\le m}$, where
$$w_i=\prod_{\substack{j\neq i,  1\leq j\leq m  }} (t_1-u_j)
\sum_{l=0}^{m-1} (t_1t_2)^l (1+s_1).$$ Since $\prod_{i=1}^{m}
(t_1-\xi^i)=0$,
$$w_i=
 \prod_{\substack{j\neq i\\
1\leq j\leq m  }} (t_1-u_j)\prod_{1\le j\le m-1}
(u_it_2-u_j)(1+s_1).$$  By rescaling the above elements,  $\{\bold
v_i\mid 1\le i\le m\}$ is a basis of $\ind_{\mathbb Z_m\wr
B_1}^{\mathbb Z_m\wr \mathfrak S_2} 1$, where
$$\bold v_i=\prod_{j\neq i} (t_1-u_j)
\prod_{j\neq m-i} (t_2-u_j)(1+s_1).$$ We have:
\begin{itemize}\item  $F\bold v_m$ is an $FW_{m, 2}$-module with $s_1\bold v_m =t_1\bold
v_m=\bold v_m$. By \cite[2.1]{DR1}, $F\bold v_m\cong S^{\eta_m}$.

\item Suppose $2\nmid m$. If $\frac{m+1} 2\le i\le m-1 $, then $\xi^{i}\neq \xi^{m-i}$.
The subspace $F \bold v_i\oplus F\bold v_{m-i}$ is an $FW_{m,
2}$-module such that $t_1\bold v_j =u_j \bold v_j$ for $j=i, m-i$,
and $s_1\bold v_{i} = \bold v_{m-i}$. Therefore, $F \bold
v_i\oplus F\bold v_{m-i}\cong S^{\eta_i}$, $\frac{m+1} 2\le i\le
m-1 $.

\item Suppose $2\mid m$. If $\frac m2<i\le  m-1$, then $F
\bold v_i\oplus F\bold v_{m-i}$ is an $FW_{m, 2}$-module such that
$ t_1\bold v_j =u_j \bold v_j$ for $j=i, m-i$, and $s_1\bold v_{i}
= \bold v_{m-i}$. Therefore, $F \bold v_i\oplus F\bold
v_{m-i}\cong S^{\eta_i}$, $\frac m2<i\le  m-1$.

\item Suppose  $i=\frac m 2$.  Then $F\bold v_{i}$ is an $FW_{m, 2}$-module such
that $s_1\bold v_{i} =\bold v_{i}$ and $t_1\bold v_{i} =u_{i}
\bold v_{i}$. Therefore, $F\bold v_{i} \cong S^{\eta_i}$.
\end{itemize}
Consequently,   $\ind_{\mathbb Z_m\wr B_1}^{\mathbb Z_m\wr \mathfrak
S_2} 1\cong \bigoplus_{\boldsymbol \eta\in \wp_m(2)} S^{\boldsymbol
\eta}$  no matter whether $m$ is even or odd.
\end{proof}

\begin{Remark} Proposition~\ref{eta} is a special case
of \cite[(4.4)]{RY}. The decomposition given there involves certain
$m$-partitions $\eta$. In fact,  we have to put more restrictions on
$\eta$. The reason is that $\sum_{l=0}^{m-1} \tilde t_i^l w
e_{\boldsymbol a}$ may be equal to zero for general $\boldsymbol a$
(Here, we keep the notation in \cite{RY}). Therefore, the first
equality in \cite[(4.3)]{RY} is not true in general. If we denote by
$c_{\eta}$ the multiplicity of $S^{\eta}$ in $\ind_{\mathbb Z_m\wr
B_k}^{\mathbb Z_m \wr \mathfrak S_{2k}} 1$, \cite[(4.1), 6.2]{RY}
are still true although we do not know the explicit description of
$c_{\eta}$. Proposition~\ref{eta} gives us the explicit information
for $\eta$ and $c_{\eta}$ when $k=1$.
\end{Remark}

In the remaining  part of this section,  we recall the result in
\cite{RY}, which says that $\cba{n}$ is a cellular algebra in the
sense of \cite{GL}. We also prove Theorem~\ref{admissible}, which
will play the key role in the proof of Theorem~\ref{RANK}.

 Recall  that a \textsf {dotted Brauer diagram} $D$ with $k$ horizontal
 arcs is determined by a
pair of labelled $(n, k)$-parenthesis diagrams $\alpha, \beta$ and
$w\in W_{m, n-2k}$, and vice versa \cite{RY}. In this situation, we
write $D=\alpha\otimes w\otimes \beta$
 if
 \begin{itemize} \item $\alpha$ (resp. $\beta$ ) is the top (resp.
bottom) row of $D$.
\item $w$ corresponds to the dotted  Brauer diagram (or braid diagram)
which is obtained from $D$ by removing the horizontal arcs at top
and bottom rows of $D$.
\end{itemize}

We denote by  $P(n, k)$  the set of all labelled $(n,
k)$-parenthesis diagrams.

Suppose $\lambda\in \Lambda_m^+(n)$. A $\lambda$-tableau is a
bijection $\t=(\t_1, \cdots, \t_{m-1}, \t_m): (Y(\lambda^{(1)}),
\cdots Y(\lambda^{(m-1)}), Y(\lambda^{(m)}))\rightarrow \{1, 2,
\cdots, n\}$. If the entries in each $\t_i$, $1\le i\le m$ increase
from left to right in each row and from top to bottom in each
column, then $\t$ is called a standard $\lambda$-tableau. Let
$T^s(\lambda)$ be the set of all standard $\lambda$-tableaux. Let
$\{y_{\s\t}^\lambda\mid \lambda\in \Lambda_m^+(n), \s, \t\in
T^s(\lambda)\}$ be the Murphy basis for $FW_{m,
n}$\cite[2.8]{DR1}.
Define
\begin{equation}\label{cell}
C_{(\alpha, \s), (\beta, \t)}^{(k, \lambda)}=\alpha\otimes
y_{\s\t}^\lambda \otimes \beta,\quad \alpha, \beta\in P(n, k), \s,
\t\in T^\s(\lambda)
\end{equation}

Recall that $R$ is a commutative ring containing the identity  $1$
and $\delta_1, \cdots, \delta_m$.
\begin{Theorem}\cite[5.11]{RY}\label{cel}
 Suppose  $R$ contains $ u_1,
\cdots, u_m$  such that $x^m-1=(x-u_1)(x-u_2)\cdots (x-u_m)$. Let
$\La=\{(f, \lambda)\mid 0\le f\le \lfloor n/2\rfloor,
\lambda\in\La_m^+(n-2f)\}$. Then
$$\{C_{(\alpha, \s), (\beta, \t)}^{(k, \lambda)}\mid \alpha, \beta\in P(n,
k), \s, \t\in T^s(\lambda), (k, \lambda)\in \Lambda\}$$ is a
cellular basis of $\cba{n}$. The $R$-linear anti-involution defined
on $\cba{n}$ is that  defined in Lemma~\ref{quasi}.
\end{Theorem}

Following \cite[2.1]{GL}, we have the \textsf {cell modules} for
$\cba{n}$ with respect to the cellular basis provided in
Theorem~\ref{cel}. Let $\Delta(k, \lambda)$ be the cell module for
$\cba{n}$ with respect to $(k, \lambda)\in \Lambda$. Let
$\Delta(\lambda)$ be the cell module for $FW_{m, n}$ with respect to
the cellular basis $\{y_{\s\t}^\lambda\mid \lambda\in
\Lambda_m^+(n), \s, \t\in T^s(\lambda)\}$.

It has been proved in \cite[2.7]{DR1} that  $\Delta(\lambda)\cong
S^{\lambda'}$, where $\lambda'$ is the dual partition of $\lambda$.
 By \cite[2.1]{GL},  $\Delta(k, \lambda)$ is spanned by
$\alpha\otimes \bold v_j\otimes\alpha_0 \mod \cba{n}^{>(k,
\lambda)}$, where $\bold v_j$ ranges over the basis elements of
$S^{\lambda'}$.

Suppose $\lambda\in \Lambda_m^+(n)$ and $\mu\in \Lambda_m^+(n-1)$.
If there is a pair $(i, j)$ such that
$\lambda^{(j)}_i=\mu^{(j)}_i+1$ and $\lambda_l^{(k)}=\mu_l^{(k)}$
for any $(k, l)\neq (i, j)$, then we write $\mu\rightarrow
\lambda$ and say that  $\mu$ is  obtained from $\lambda$ by
removing a box.  In this situation, we also say that $\lambda$ can
be obtained from $\mu$ by adding a box.

\begin{Theorem} \label{admissible} Let $\cba{n}$ be a cyclotomic
Brauer algebra over $F$. If $\mu\in \Lambda_m^+(n-2)$ and
 $\lambda\in \Lambda_m^+(n)$, then either  $[\Delta(1, \mu'): \Delta(\lambda')]=0$ or
$[\Delta (1, \mu'): \Delta(\lambda')]=1$. Furthermore, $[\Delta(1,
\mu'): \Delta(\lambda')]=1$ if one of the following conditions holds
true.
\begin{itemize}
\item [(1)] $\lambda^{(j)}=\mu^{(j)}$, $j\neq m$ and two  boxes in the
skew Young diagram $Y(\lambda^{(m)}/\mu^{(m)})$ are not in the same
column.

\item [(2)] Suppose that $m$ is odd. There is an $i$ with $\frac {m+1}2\le i\le m-1$ such that
$\mu^{(i)}\rightarrow \lambda^{(i)}$ and $\mu^{(m-i)}\rightarrow
\lambda^{(m-i)}$, and $\lambda^{(j)}=\mu^{(j)}$ for $j\neq i, m-i$.

\item [(3)] Suppose that $m$ is even. There is an $i$ with $\frac {m}2<  i\le m-1$ such that
$\mu^{(i)}\rightarrow \lambda^{(i)}$ and $\mu^{(m-i)}\rightarrow
\lambda^{(m-i)}$, and $\lambda^{(j)}=\mu^{(j)}$ for $j\neq i, m-i$.

\item [(4)] Suppose that $m$ is even. $\lambda^{(j)}=\mu^{(j)}$, $j\neq m/2$
and two  boxes in the skew Young diagram $Y(\lambda^{(\frac
m2)}/\mu^{(\frac m2)})$ are not in the same column.
\end{itemize}
\end{Theorem}

\begin{proof} Since we are assuming that $F$ is a splitting field of $x^m-1$
such that $e \nmid m\cdot n!$, both $FW_{m, k}$ and $F\mathfrak S_k$
are (split) semisimple for $k\le n$.

For  $\lambda\in\La^+(n_1+n_2) , \mu\in \Lambda^+(n_1), \eta\in
\Lambda^+(n_2)$, let $L_{\eta, \mu}^\lambda$ be the corresponding
Littlewood-Richardson coefficient for symmetric groups. If
$\lambda\in \La_m^+(n)$, $\mu\in \Lambda_m^+(n-2)$ and $\eta\in
\La_m^+(2)$, the Littlewood-Richardson coefficient $L_{\eta,
\mu}^\lambda$ for complex reflection groups is $\prod_{i=1}^m
L_{\eta^{(i)}, \mu^{(i)}}^{\lambda^{(i)}}$ \cite[\S 4]{RY}.  Let
$c_\eta$ be the multiplicity of $S^{\eta}$ in $\ind_{\mathbb Z_m
\wr B_1}^{W_{m, 2}} 1$. By \cite[6.2]{RY}, $[\Delta (1, \mu'):
\Delta(\lambda')]=m_{\mu, \lambda}$, where $m_{\mu,
\lambda}=\sum_{\eta\in \wp_m(2)} c_\eta L_{\eta, \mu}^{\lambda}$.
Note that $L_{\eta, \mu}^\lambda\neq 0$ if and only if
$L_{\eta^{(i)}, \mu^{(i)}}^{\lambda^{(i)}}\neq 0$ for $1\le i\le m$.
Consequently, if  $L_{\eta, \mu}^\lambda\neq 0$ for $\eta\in
\wp_m(2)$, then  there is a unique  $\eta\in \wp_m(2)$. Suppose
$\lambda$ and $\mu$ are partitions. It is known that two boxes in
the skew Young diagram $Y(\lambda/\mu)$ are not in the same column
if $ L_{(2), \mu}^\lambda\neq 0$ (see, e.g. \cite[\S3]{DWH}). In
this situation, $ L_{(2), \mu}^\lambda=1$, and $\lambda\supset \mu$.
We use classical branching rule for symmetric groups if either
$\eta\not\in \{\eta_m, \eta_{\frac m2}\}$, $2\mid m$ or $\eta\neq
\eta_m$, $2\nmid  m$. In any case, we have $m_{\mu, \lambda}=1$, if
one of conditions in (1)-(4) holds true.
\end{proof}

\begin{Defn}\label{adm} Suppose $\mu\in \Lambda_m^+(n-2)$ and
$\lambda\in \Lambda_m^+(n)$. $\lambda$  is called
\text{$\mu$-admissible} if one of the conditions in
Theorem~\ref{admissible}(1)-(4) holds true. Let $\mathscr A(\mu)$ be
the set of all $\mu$-admissible $m$-partitions.
\end{Defn}

\section{Zero divisors of certain discriminants}

In this section, we assume  $\delta_i\in F$ for $1\le i\le m$, where
$F$ is a splitting field of $x^m-1$ and $e\nmid m\cdot n!$. The main
purpose of this section is to prove Theorem~\ref{key}, which will
give all zero divisors of the discriminants of the Gram matrices
$\G_{1, \mu'}$ with respect to the cell modules $\Delta(1, \mu')$,
$\mu\in \La_m^+(n-2)$.

Recall that $P(n, k)$ is the set of labelled parenthesis Brauer
diagrams with k horizontal arcs. In what follows, we assume
$\alpha_0=top(e_{n-1})\in P(n, 1)$, the top row of $e_{n-1}$. Define
$M_1$ and $M_2$ by setting
\begin{itemize}
\item  $M_1=\set {\alpha\otimes w\otimes \alpha_0\mid \alpha\in P(n,
1), w\in W_{m, n-2}}.$

\item  $M_2=\set {\alpha\otimes w\otimes \beta\mid \alpha, \beta\in
P(n, k), w\in W_{m, n-2k}, 2\le k\le \lfloor \frac n2\rfloor}$.
\end{itemize}

We consider the quotient $F$-subspace $V=V_1/V_2$, where $V_1$
(resp. $V_2$)  is spanned by $M_1\cup M_2$ (resp. $M_2$). For
convenience, we use $\alpha\otimes w\otimes \alpha_0$ instead of
$\alpha\otimes w\otimes \alpha_0+V_2$.

Recall that any dotted  Brauer diagram can be written as $
\alpha\otimes w\otimes \beta $ where $\alpha, \beta\in P(n, k)$ and
$w\in W_{m, n-2k}$. Let $\tilde \alpha\in P(n, k)$ be  such that
\begin{enumerate}\item $\alpha$ and $\tilde \alpha$ have the same
horizontal arcs.
\item There are $m-i$ dots on a horizontal arc  in $\tilde \alpha$
if and only if there are $i$ dots on the corresponding horizontal
arc in $\alpha$.\end{enumerate}

 Define an $R$-linear isomorphism  $\iota: \cba{n}\longrightarrow\cba{n}$ by declaring that
 \begin{equation}\label{iota}\iota (\alpha\otimes w\otimes \beta)=\tilde
\beta\otimes w^{-1} \otimes \tilde \alpha.\end{equation}
 We remark that $\iota$
is not an algebraic (anti-)homomorphism since $\iota (e_i t_i^k
e_i)=\delta_k e_i\neq \delta_{m-k} e_i$ in general. However, by
straightforward computation, we have
\begin{equation}\label{gaction} \iota (w (\alpha \otimes w_1\otimes
\beta))= \iota(\alpha\otimes w_1\otimes \beta) w^{-1},
\end{equation}
 for any $\alpha, \beta\in P(n, k), w\in W_{m,
n}, w_1\in W_{m, n-2k}$.

Following \cite{HW1}, we have the following definition.

\begin{Defn} Suppose  $\alpha_i\otimes w\otimes \alpha_0\in V$
for $i=1,2$. Let $\langle \alpha_1\otimes w_1\otimes \alpha_0,
\alpha_2 \otimes w_2\otimes \alpha_0\rangle$ be the coefficient of
$e_{n-1}$ in the expression of $\iota (\alpha_1\otimes w_1\otimes
\alpha_0) \cdot ( \alpha_2 \otimes w_2\otimes \alpha_0)$, where
$\iota$ is defined in (\ref{iota}).
 Let $\G_{m, n}(\boldsymbol\delta)$ be the $f\times f$-matrix
with $f=\dim V$ such that the entry in $(\alpha_1\otimes w_1\otimes
\alpha_0)$-th row, $(\alpha_2 \otimes w_2\otimes \alpha_0)$-th
column is $\langle \alpha_1\otimes w_1\otimes \alpha_0, \alpha_2
\otimes w_2\otimes \alpha_0\rangle$. \end{Defn}

If either $h_1\in M_2$ or $h_2\in M_2$, then   $\iota (h_1) h_2\in
V_2$. Since  $e_{n-1}\not\in V_2$, $\langle h_1, h_2\rangle=0$.
Hence, $\langle\ \ , \ \ \rangle: V\times V\longrightarrow F$  is a
well-defined $F$-bilinear form on $V$.

The following lemma can be verified easily.

\begin{Lemma} $\G_{m, n}(\boldsymbol\delta)=(g_{ij})$ is an
$f\times f$ matrix such that $g_{ii}=\delta_0$, $1\le i\le f$ and
$g_{ij}\in \set { 0, 1, \delta_1, \cdots, \delta_{m-1}}$ if $i\neq
j$.
\end{Lemma}

\begin{Lemma} $\G_{m,
n}(\boldsymbol\delta):  V\rightarrow  V$ is a  left $FW_{m,
n}$-homomorphism and a right $FW_{m,n-2}$ homomorphism
\end{Lemma}

\begin{proof} We consider  $\G_{m,
n}(\boldsymbol\delta)$ as the $F$-linear endomorphism on $V$ such
that
$$\G_{m,
n}(\boldsymbol\delta) (\alpha_1\otimes w_1\otimes
\alpha_0)=\sum_{\alpha\in P(n, 1), w\in W_{m, n-2} }\langle
\alpha\otimes w\otimes \alpha_0, \alpha_1 \otimes w_1\otimes
\alpha_0\rangle \alpha\otimes w\otimes \alpha_0.$$ By
(\ref{gaction}),
$$\langle w
(\alpha_1\otimes w_1\otimes \alpha_0), w (\alpha_2 \otimes
w_2\otimes \alpha_0)\rangle=\langle \alpha_1\otimes w_1\otimes
\alpha_0, \alpha_2 \otimes w_2\otimes \alpha_0\rangle.$$ In other
words, $\G_{m, n}(\boldsymbol\delta): V\rightarrow V$ is a left
$FW_{m, n}$-homomorphism.

On the other hand, since $(\alpha_1\otimes w_1\otimes \alpha_0) y=
\alpha_1\otimes w_1y\otimes \alpha_0$ for any $y\in W_{m, n-2}$,
$e_{n-1}$ appears in $y^{-1}(\alpha_0\otimes w_1\otimes
\alpha_0)y$ with non-zero coefficient if and only if $w_1=1$.
Therefore,
$$\langle
(\alpha_1\otimes w_1\otimes \alpha_0) y,  (\alpha_2 \otimes
w_2\otimes \alpha_0) y\rangle=\langle \alpha_1\otimes w_1\otimes
\alpha_0, \alpha_2 \otimes w_2\otimes \alpha_0\rangle.$$
Consequently, $\G_{m, n}(\boldsymbol\delta):  V\rightarrow  V$ is
a right $FW_{m, n-2}$-homomorphism.
\end{proof}

Since we are assuming that  $F$ is a splitting field of $x^m-1$ and
 $e\nmid m\cdot n!$, $FW_{m, k}$ is (split) semisimple for any $k$,  $1\le
k\le n$. Assume that $\lambda\in \Lambda_m^+(k)$. The classical
Specht module $S^{\lambda}$ is a direct summand of $FW_{m, k}$.
Consequently, $\Delta(1, \lambda')$ can be realized as a submodule
of $V$, which is spanned by $\alpha\otimes v_j\otimes \alpha_0 (\mod
V_2)$, where $v_j$ ranges over
 the basis elements of  $S^{\lambda}$. Note that
 $\G_{m, n}(\boldsymbol \delta)$ is a right $FW_{m,
n-2}$-module. For any $\lambda\in \La_m^+(n-2)$, the restriction of
$\G_{m, n}(\boldsymbol \delta)$ on $\Delta(1, \lambda')$ induces a
linear endomorphism on $\Delta(1, \lambda')$.

\begin{Defn}
For  $\mu\in \Lambda_m^+(n-2)$, define
$g_\mu=\prod_{\lambda\in\mathscr A(\mu)} g_{\lambda, \mu}$, where
\begin{equation} g_{\lambda, \mu}=(\bar \delta_0-m+m\sum_{p\in
Y(\lambda/\mu)} c(p)) \prod_{i=1}^{m-1} (\bar \delta_i+m\sum_{p\in
Y(\lambda/\mu)} c(p)).\end{equation}
\end{Defn}

It follows from \cite[2.3]{GL} that there is an invariant symmetric
bilinear form defined on each cell module $\Delta(k, \lambda)$. Via
such a bilinear form, one can define a Gram matrix $\G_{k,
\lambda}$. Let $\det \G_{k, \lambda}$ be the determinant of $\G_{k,
\lambda}$. The following result follows from \cite[3.8]{GL} and
Theorem~\ref{cel}, immediately.

\begin{Lemma}  $\cba{n}$
is (split)  semisimple over $F$  if and only if $\det \G_{k,
\lambda}\neq 0$ for all $(k, \lambda)\in \Lambda$.
\end{Lemma}

In general,  it is difficult to compute $\det \G_{k, \lambda}$.
Assume  $\delta_i\neq 0$ for some $1\le i\le m$. The following
result describes all the zero divisors of $\det \G_{1, \lambda}$,
$\lambda\in \Lambda_m^+(n-2)$ . Fortunately, it completely
determines $\cba{n}$ being (split) semisimple.

\begin{Theorem}\label{key}  Suppose  $\delta_i\neq 0$ for some $i, 1\le i\le m$. $\det \G_{1,
\mu'}= 0$ if and only if $g_{\mu}=0$.
\end{Theorem}

\begin{proof}
($\Rightarrow$)  If $\det \G_{1, \mu'}=0$, then we can find an
irreducible
 $\cba{n}$-module $M\subset
\Rad \G_{1, \mu'}$, where $\Rad \G_{1, \mu'}=\set{v\in \Delta(1,
\mu')\mid \G_{1, \mu'} (v)=0}$. It follows from \cite[2.6, 3.4]{GL}
that any irreducible module of a cellular algebra must be the simple
head of a cell module, say $\Delta(k, \lambda')$. Hence,  there is
a non-zero homomorphism from $\Delta(k,\lambda')$ to $\Delta(1,
\mu')$ with $(k, \lambda')< (1, \mu')$.  Therefore, either $k=1$ or
$k=0$.

Assume that $k=1$. We use \cite[7.4]{RY}\footnote{\cite[7.4]{RY} is
for $\cba{n}$ over the complex field. However, it is still true if
we use $F$, a splitting field of $x^m-1$ with $e\nmid m\cdot n!$,
instead of $\mathbb C$. See \cite[8.8]{RY}.} to get a non-zero
homomorphism from $\Delta(0, \lambda')$ to $\Delta(0, \mu')$. Notice
that, as $FW_{m, n}$-modules, $\Delta(0, \lambda')\cong
S^{\lambda}$. Since  $FW_{m, n}$ is (split) semisimple, we have
$\lambda=\mu$, a contradiction since $(1, \lambda')<(1, \mu')$.

If $k=0$, then there is a non-zero $\cba{n}$-homomorphism from
$\Delta(0, \lambda')$ to $\Delta(1, \mu')$, forcing $\lambda\in
\mathscr A(\mu)$. By \cite[8.6, 8.8]{RY},  $g_{\lambda, \mu}=0$. We
have $g_\mu=0$ as required.

$(\Leftarrow)$ Suppose $g_\mu=0$. Then there is a $\lambda\in
\mathscr A(\mu)$ such that $g_{\lambda, \mu}=0$. Since $\lambda\in
\mathscr A(\mu)$, by Theorem~\ref{admissible}, $[\Delta(1, \mu'):
S^\lambda]=1$. Hence,   there is a unique  $FW_{m, n}$-submodule $M$
of $\Delta(1, \mu')$ which is isomorphic to $S^\lambda$. Recall that
$\G_{m, n}(\boldsymbol \delta)\mid_{\Delta(1, \mu')}$ is a linear
endomorphism on $\Delta(1, \mu')$. For  simplicity, we use $\G_{m,
n}(\boldsymbol \delta)$ instead of  $\G_{m, n}(\boldsymbol
\delta)\mid_{\Delta(1, \mu')}$ if there is no confusion.

Since  $\G_{m, n}(\boldsymbol \delta)$ is an $FW_{m,
n}$-homomorphism, and $[\Delta(1, \mu'): S^\lambda]=1$, $\G_{m,
n}(\boldsymbol \delta) (M)\subset M$. By Schur's Lemma, $\G_{m,
n}(\boldsymbol \delta)\mid_M =f(\boldsymbol \delta)I$, where $I$ is
 $\dim M\times \dim M$ identity matrix and  $f(\boldsymbol
\delta):=f(\delta_0, \delta_1, \cdots, \delta_{m-1})$ is a
polynomial in $\delta_i, 0\le i\le m-1$.

Take a basis of $M$ and extend it to get a basis of $ V$ via the
elements $\alpha\otimes w\otimes \alpha_0$. Then $\G_{m,
n}(\boldsymbol \delta)$ is conjugate to
$\begin{pmatrix} f(\boldsymbol \delta) I & 0\\
\ast & B\\
\end{pmatrix},$
where  any entry in the diagonal of $B$ is $\delta_0$,  and the
term of the entry of $B$ elsewhere does not contains $\delta_0$.
 Since the  degree of $\delta_0$ in $\det \G_{m, n}(\boldsymbol
\delta)$ is $\dim V$ (see Lemma~3.2), the degree of $\delta_0$ in
$f(\boldsymbol \delta)$ must be  $1$. In particular,  $f(\boldsymbol
\delta)$ is not a constant number. Take the parameters $\delta_0,
\delta_1, \cdots, \delta_{m-1}$ such that $f(\boldsymbol \delta)=0$.
Then $\G_{m, n}(\boldsymbol \delta) \mid_{M}=0$.

We claim $e_{n-1} v=0$ for any $v\in M$. Write $v=\sum_{\alpha^s, w}
a_{\alpha^s, w} \alpha^s\otimes w\otimes \alpha_0$,  where there are
$s$ dots at the left endpoint of the unique arc in $\alpha^s$. We
divide $P(n, 1)$  into three disjoint subsets $P_1, P_2, P_3$ as
follows. Recall that  a point in $\alpha^s$ is called a fixed  point
if it is an endpoint of a horizontal  arc of $\alpha^s$. Otherwise,
it is called a free point.

\begin{itemize}
\item  $P_1$ consists of all $\alpha^s\in P(n, 1)$ such that $(n-1, n)$ is
a unique arc of $\alpha^s$.  Then $e_{n-1} (\alpha^s\otimes
w\otimes \alpha_0)=\delta_s \alpha_0\otimes w\otimes \alpha_0$.

\item $P_2$ consists of all $\alpha^s\in P (n, 1)$ such
that both $n-1$ and $n$ are free points in $\alpha$. Then $e_{n-1}
\alpha^s\otimes w\otimes \alpha_0=0$.

\item  $P_3$ consists of all $\alpha^s\in P(n, 1)$ such that
either  $n-1$ or $ n$ is a fixed point. Let $i$  be the left
endpoint of the unique arc in $\alpha^s$. By assumption, there are
$s$ dots at the endpoint $i$. We define $w_{\alpha^s} \in
\mathfrak S_{n-2}$ by setting
$$w_{\alpha^s}=\begin{pmatrix} i & i+1 & i+2 & \cdots & n-3 & n-2\\
                        n-2 & i & i+1 & \cdots & n-4 & n-3\\
                        \end{pmatrix}
                        $$
                        \end{itemize}
Define $y_{\alpha^s}:=t_i^s w_{\alpha^s}$. Then $e_{n-1}\cdot
(\alpha^s\otimes 1\otimes \alpha_0)= \alpha_0\otimes
y_{\alpha^s}\otimes \alpha_0$. Therefore, the coefficient of
$\alpha_0 \otimes w_1\otimes \alpha_0$ in $e_{n-1} v$  is
$\sum_{\alpha^s\in P_3} a_{\alpha^s, y_{\alpha_s}^{-1} w_1}
+\sum_{s=0}^{m-1} \delta_\s a_{\alpha_0^s, w_1}$.

On the other hand, by direct computation, the coefficient of
$\alpha_0\otimes w_1\otimes \alpha_0$ in $\G_{m, n}  (\boldsymbol
\delta) v$ is $ \sum_{\substack {\alpha^s\in P(n, 1),\\ w\in W_{m,
n-2}}} a_{\alpha^s, w} \langle \alpha_0\otimes w_1\otimes
\alpha_0, \alpha^s \otimes w\otimes \alpha_0 \rangle $. We have $$
\langle \alpha_0\otimes w_1\otimes \alpha_0, \alpha^s \otimes
w\otimes
\alpha_0 \rangle=\begin{cases} \delta_s, &\text{if $\alpha^s\in P_1, w=w_1$,}\\
                                0, &\text{if $\alpha^s\in P_1, w\neq
                                w_1$,}\\
                          0,         &\text{if $\alpha^s\in P_2$,}\\
                          1, &\text{if $\alpha^s\in P_3$ and
                          $w=y_{\alpha^s}^{-1} w_1$.}\\
                          0, &\text{if $\alpha^s\in P_3$ and
                          $w\neq y_{\alpha^s}^{-1} w_1$.}
                          \end{cases}$$
                          Since $\G_{m, n} (\boldsymbol
\delta)v=0$, the coefficient
                          of $\alpha_0\otimes w_1\otimes
                          \alpha_0$ in $\G_{m, n}(\boldsymbol
\delta)v$ is zero.
                          Therefore,
                           $${\tiny
\sum_{\substack {\alpha^s\in P(n, 1),\\ w\in W_{m, n-2}}}
a_{\alpha^s, w} \langle \alpha_0\otimes w_1\otimes \alpha_0,
\alpha^s \otimes w\otimes \alpha_0 \rangle = \sum_{\alpha^s\in
P_3} a_{\alpha^s, y_{\alpha^s}^{-1} w_1} +\sum_{s=0}^{m-1}
\delta_\s a_{\alpha_0^s, w_1}=0,}
$$
forcing the  coefficient of $\alpha_0\otimes w_1\otimes \alpha_0$
in $e_{n-1} v$ to be zero for all $w_1\in W_{m, n-2}$. This
completes the proof of the claim.

Therefore,   as $\cba{n}$-module, $M\cong \Delta(0, \lambda')$. We
obtain  a non-zero $\cba{n}$-homomorphism from $\Delta(0, \lambda')$
to $\Delta(1, \mu')$. In particular, $\det G_{1, \mu'}=0$. By
\cite[8.6, 8.8]{RY}\footnote{under our assumption, the group algebra
$FW_{m, n}$ is (split) semisimple. Since the proof of \cite[8.6,
8.8]{RY} depends only  on the fact that $\mathbb C W_{m, n}$ is
(split) semisimple, we can apply these results here.} the parameters
$\delta_i$'s  must satisfy the equation $g_{\lambda, \mu}=0$, the
condition we have assumed.
\end{proof}

\section{Proof of Theorem~\ref{RANK} }

In this section, we prove Theorem~\ref{RANK}, the main result of
this paper. Unless otherwise stated, we assume that $F$ is a
splitting field of $x^m-1$, which contains $\delta_i, 1\le i\le m$.
Assume  $m>1$.

\begin{Prop} \label{81} Suppose $n\ge 2$. If $0\neq \delta_i\in F$ for some $i$, $1\le
i\le m$ , then $\cba{n}$   is (split) semisimple if and only if  $ e
\nmid m \cdot n!$ and $\det \G_{1, \lambda}\neq 0$ for any
$\lambda\in \Lambda_m^+(k-2)$, $2\le k\le n$.
\end{Prop}

\begin{proof}
$(\Leftarrow)$ Suppose that $\cba{n}$ is not (split) semisimple.
There is a $(k, \lambda)\in \Lambda$  such that $\det \G_{k,
\mu}=0$.
 Since $FW_{m, n}$ is (split)
semisimple, $k\neq 0$.

Take an irreducible submodule $M\subset \Rad \Delta(k, \mu)$. By
\cite[2.6, 3.4]{GL}, $M$ must be isomorphic to the simple head of a
cell module, say $\Delta(l, \lambda)$, such that  $(l, \lambda)<(k,
\mu)$. Furthermore, it results in a non-trivial homomorphism from
$\Delta(l, \lambda)$ to $\Delta(k, \mu)$.

If $l=k$, we use  \cite[7.4]{RY} to  get $\Delta(0, \la)\cong
\Delta(0, \mu)$. As $FW_{m, n-2k}$-modules, $\Delta(0, \lambda)\cong
S^{\lambda'}$. Since $FW_{m, n-2k}$ is (split) semisimple,
$\la=\mu$, which  contradicts  $(l, \lambda)<(k, \mu)$.

 Suppose $l<k$. By
\cite[7.4, 7.7]{RY}, there is a non-trivial homomorphism from
$\Delta(0, \tilde \la)$ to $\Delta(1, \tilde \mu)$ for some
$\tilde \mu\in \La_m^+(p-2)$ with $p\le n$. By assumption,
 $\det \G_{1, \tilde\mu}\neq 0$. Hence,
 $\Delta(1, \tilde \mu)=D^{(1, \tilde \mu)}\cong \Delta(0, \tilde
 \la)$. By \cite[3.4]{GL},   $(0, \tilde \la)= (1, \tilde \mu)$,
 a contradiction.

$(\Rightarrow)$ If $\cba{n}$ is (split) semisimple, then
\cite[3.8]{GL} implies that $\det \G_{k, \lambda}\neq 0$ for all
$0\le k\le \lfloor\frac n2\rfloor$. Therefore, $FW_{m, n}$ is
(split) semisimple, forcing $e\nmid m\cdot n!$.

Suppose $\det \G_{1, \mu'}=0$ for some $\mu' \in \Lambda_m^+(k-2)$.
Then $k<n$. By Theorem~\ref{key}, there is a $\mu$-admissible
$m$-partition $\lambda$ such that $g_{\lambda, \mu}=0$.
Equivalently, there is a non-zero $\cba{n}$-homomorphism from
$\Delta(0, \lambda')$ to $\Delta(1, \mu')$.

Since we are assuming that  $m\ge 2$,  we can find an $i, 1\le i\le
m$, such that $\lambda^{(i)}=\mu^{(i)}$. We can add $l$ boxes to
$\lambda^{(i)}$ so as to  get another partition $\tilde
\lambda^{(i)}=\tilde\mu^{(i)}$. In this situation,
$g_{\tilde\lambda, \tilde\mu}=g_{\lambda, \mu}$, where
$\tilde\lambda$ (resp. $\tilde\mu$ ) can be obtained from $\lambda$
(resp. $\mu$) by using $\tilde \lambda^{(i)}$ instead of
$\lambda^{(i)}$ (resp. $\mu^{(i)}$). By definition, $\tilde\la\in
\mathscr A(\tilde\mu)$. If we take $l$ such that $|\lambda|+l=n$,
then $\Delta(0, {\tilde\lambda}')$ and $\Delta(1,{\tilde\mu}')$ are
$\cba{n}$-modules. By Theorem~\ref{key}, $\det \G_{1,
{\tilde\mu}'}=0$. However, since $\cba{n}$ is (split) semisimple,
$\det \G_{1, {\tilde\mu}'}\neq 0$, a contradiction.
\end{proof}

\begin{Cor} \label{problem} Let $\cba{n}$ be a cyclotomic Brauer
algebra over $F$, where $F$ contains a non-zero  $\delta_i$ for some
$i, 1\le i\le m$. $\cba{n}$ is (split) semisimple if only if $\det
\G_{k, \lambda}\neq 0$ for all
 $\lambda\in \Lambda_m^+(n-2k)$, and $k=0, 1 $.\end{Cor}

\begin{proof} Suppose $\cba{n}$ is (split) semisimple. It follows from
\cite[3.8]{GL} that $\det \G_{k, \lambda}\neq 0$ for all $0\le
k\le \lfloor\frac n2\rfloor$. In particular, $\det \G_{k,
\lambda}\neq 0$ with $k=0, 1 $ and $\lambda\in \Lambda_m^+(n-2k)$.

Conversely, if $\det \G_{0, \lambda}\neq 0$ for all $\lambda\in
\Lambda_m^+(n)$, then $FW_{m, n}$ is (split) semisimple. Suppose
that $\cba{n}$ is not (split) semisimple.  By Proposition~\ref{81},
there is a $\mu\in \Lambda_m^+(k-2)$ with $k<n$ such that $\det
\G_{1, \mu}= 0$. From the proof of Proposition~\ref{81},   we can
find a $\tilde \mu\in \La_m^+(n-2)$ such that  $\det \G_{1,
\tilde\mu}=0$. This contradicts our assumption.
\end{proof}

Corollary~\ref{problem} has been stated  as a question in
\cite[p220]{RY}  . We remark that corollary~\ref{problem}  is not
true if $m=1$. In fact, the first author has  proved that a Brauer
algebra  is  (split) semisimple over $F$  if and  only if $e\nmid
n!$ and $\det G_{1, \lambda}\neq 0$ for all  $\lambda\in
\La^+(k-2)$, $2\le k\le n$. By \cite[3.3-3.4]{DWH},
Corollary~\ref{problem} is not true if  $m=1$.

\begin{Defn} Suppose  that $m, n\in \mathbb N$ with $ n\ge 2$. For
$m\ge 2$, define
 $\rho_{m, n}=\{ma\mid a\in \tilde\rho_{m, n}\}$, where $$ \tilde\rho
_{m, n}=\{k\in \mathbb Z\mid k=\sum_{p\in Y(\la/\mu)} c(p)\mid
\mu\in \La_m^+(n-2), \lambda\in \mathscr A(\mu)\}.$$
 If $m=1$, we
define $$ \tilde{\rho}_{m, n}=\{r\in \mathbb Z\mid r=\sum_{p\in
Y(\la/\mu)} c(p)\mid \mu\in \La^+(k-2), \lambda\in \La^+(k), 2\le
k\le n\},$$ where two boxes in $Y(\lambda/ \mu)$ are not in the same
column.
\end{Defn}
At the end of this paper, we will prove $\tilde\rho_{m,
n}=\tilde{\mathbb Z}_{m, n}$. Hence,  $\rho_{m, n}=\mathbb Z_{m,
n}$.

\begin{Theorem}\label{main}
Let $\cba{n}$ be a cyclotomic Brauer algebra over $F$, where $F$
contains a non-zero  $\delta_i$ for some $i, 1\le i\le m$.
 Suppose $n\ge 2$.  $\cba{n}$ is
  (split) semisimple if and only if
\begin{itemize}\item[(1)] $e\nmid m\cdot n!$,
\item[(2)] $\epsilon_{i, 0} m-\bar\delta_i\not\in \rho_{m, n}$,
$0\le i\le m-1$, where $\epsilon_{i, 0}$ is the Kronecker
function.
\end{itemize}
\end{Theorem}
\begin{proof} The result follows from Theorem~\ref{key} and
Corollary~\ref{problem}.
\end{proof}

In the remaining  part of this section, we deal with the case
$\delta_i=0$ for all $1\le i\le m$. First, we discuss $\ccba{3}$.

 We want to compute $\det G_{1, \lambda}$ with
 $\lambda=( (1), 0, \cdots, 0)$. Note that we have assumed $u_i=\xi^i$,
 $1\le i\le m$. In this situation, $y_{\la'} w_{\la'} x_{\lambda}
=g(t_1)=\prod_{i=1}^{m-1} (t_1-\xi^i)$.
  Write
$v_1^{(0)}=top(e_1)$, $v_2^{(0)}=top (s_1e_2)$ and
$v_3^{(0)}=top(e_2)$. Let $v_i^{(k)}$ be obtained from $v_i^{(0)}$
by putting $k$ dots at the left endpoint of the unique horizontal
arc in $v_i^{(0)}$. Then $\Delta(1, \lambda)$ can be considered as a
free $F$-module with basis $\{v_i^{(k)}\otimes g(t_1)\otimes
v_3^{(0)}\mid 1\le i\le 3, 0\le k\le m-1\}$. Let
$a=\prod_{i=1}^{m-1} (1-\xi^i)$. The Gram matrix with respect to
this basis is
$$
\G_{1, \lambda}=\begin{pmatrix} 0 & A & A\\
         A &  0 & A\\
         A & A & 0\\
         \end{pmatrix},
$$
where $A=(a_{ij})$ is the $m\times m$ matrix with $a_{ij}=a$, $1\le
i, j\le m$. Since we are assuming that $m>1$, $\det \G_{1, \lambda}
= 0$. In other words, $\Rad \Delta(1, \lambda)\neq 0$. Take an
irreducible submodule $D$ of $\Rad \Delta(1, \lambda)$. Note that
any irreducible module must be the simple head of a cell module, say
$\Delta(k, \mu)$. Therefore, there is a non-trivial homomorphism
from $\Delta(k, \mu)$ to $\Delta(1, \lambda)$. By \cite[2.6]{GL},
$(k, \mu)< (1, \lambda)$. This proves the following lemma.

\begin{Lemma}\label{cba3}
Suppose  $\lambda=((1), 0, \cdots, 0)$. There is a cell module
$\Delta(k, \mu)$ of $\ccba{3}$  with $(k, \mu)< (1, \lambda)$ such
that there is a non-trivial homomorphism from $ \Delta(k, \mu)$ to
$\Delta(1, \lambda) $.
\end{Lemma}

\newcommand {\jcba}[1]{J_{m,#1}(\boldsymbol{0})}
\newcommand {\jdcba}[1]{J_{m,#1}(\boldsymbol{\delta})}
Let $\jcba{n}$ be the left ideal of $\ccba{n}$ spanned by the dotted
Brauer diagrams $D$ such that $\{n-1,n\}$ is a horizontal arc at the
bottom row  of $D$. It is clear that $\jcba{n}=\ccba{n}e_{n-1}$.

\newcommand{\iicba}[2]{I_{m,#1}^{#2}(\boldsymbol{0})}
Following \cite{RY}, let $I_{m, n}^{\ge k}$ (resp. $I_{m, n}^{>k}$)
be the vector space generated by  $(n, l)$-dotted Brauer diagrams
with  $l\ge k$  (resp. $l>k$). Let $\iicba{n}{k}= I_{m, n}^{\geq
k}/I_{m, n}^{>k}$.  Then $\iicba{n}{k}$ is a $\ccba{n}$-module. Let
$\iocba{n}{k}$  be the subspace of $\iicba{n}{k}$ generated by
$\{\alpha\otimes w \otimes \beta_0 \mid \alpha\in P(n,k), w\in W_{m,
n-2k} \}$, where $\beta_0=top (e_{n-2k+1}\cdots e_{n-3}e_{n-1})$.
Let $\ccba{n}$-mod be the category of the left $\ccba{n}$-modules.
Let
$$G: \ccba{n-2}\text{-mod}\longrightarrow
\ccba{n}\text{-mod}$$ be the tensor functor defined by declaring
that $G(M)=\jcba{n}\otimes_{\ccba{n-2}} M$, for any $\ccba{n-2}$-mod
$M$.

\begin{Prop}\label{hom} Suppose  $\lambda\in \La_m^+(n-2k)$.
\begin{itemize}
\item [(a)] The functor $G$ sends  non-zero $\ccba{n-2}$-homomorphisms to
non-zero ones.

\item[(b)] $G(\Delta(k-1,\lambda))= \Delta{(k,\lambda)}$.\end{itemize}\end{Prop}

\begin{proof} Suppose $\phi: M_1\rightarrow M_2$ is a
$\ccba{n-2}$-module homomorphism. Write $\phi_\ast=G(\phi)$. For
any $D_1\in\ccba{n}$, $D\in \jcba{n}$ and $m\in M_1$,
$$\begin{aligned} \phi_\ast (D_1 (D\otimes m))& = \phi_\ast (D_1 D\otimes
m)=(D_1 D)\otimes \phi(m)\\
&= D_1 (D\otimes \phi(m))=D_1\phi_\ast (D\otimes m)\\
\end{aligned}
$$
Therefore, $\phi_\ast$ is a $\ccba{n}$-homomorphism. For any
$\ccba{n-2}$-module $M$, define an $F$-linear map $\alpha:
\jcba{n}\otimes_{\ccba{n-2}} M\rightarrow M$ by setting
$\alpha(D\otimes m)=(e_{n-1} D)_0 m$, where  $(e_{n-1} D)_0$ is
obtained from $e_{n-1} D$ by removing the horizontal arcs $\{n-1,
n\}$ at the top and bottom rows of $e_{n-1} D$.

Suppose $D^\ast =s_{n-2} e_{n-1}\in \jcba{n}$. Then $\alpha(D^\ast
\otimes m)=m$. If $\phi\neq 0$, then there is an $m_1\in M_1$ such
that $\phi(m_1)=m_2\neq 0$. Consequently, $\alpha(D^\ast\otimes
m_2)=m_2\neq 0$. We have $\phi_\ast\neq 0$ since
$\phi_\ast(D^\ast\otimes m_1)=D^\ast\otimes m_2\neq 0$. This
completes the proof of (a).

 (b) can be proved similarly as
\cite[7.2]{RY}. We include a proof as follows. First. we claim as
$(\ccba{n}, W_{m, n-2k})$-modules
\begin{equation}
\iocba{n}{k}\cong \jcba{n}\otimes_{\ccba{n-2}}\iocba{n-2}{k-1}.
\end{equation}

For the simplification in exposition and notation,   we omit
$\ccba{n-2}$ in what follows.

 Suppose  $D_1\otimes D_2\in
\jcba{n}\otimes\iocba{n-2}{k}$. Let $e_{i, j}=\alpha\otimes 1\otimes
\alpha$, where $\alpha\in P(n, 1)$ contains a unique horizontal arc
$\{i, j\}$. Define $e_{i, j}^{s, t}=t_i^s e_{i, j}t_i^t$. We claim
that there is a dotted Brauer diagram $D_1'$ in $\iocba{n}{1}$ such
that $D_1\otimes D_2=D_1'\otimes e_{i_1,j_1}^{s_1,t_1}\cdots
e_{i_{k-1},j_{k-1}}^{s_{k-1},t_{k-1}}D_2$, where
$e_{i_l,j_l}^{s_l,t_l}\in \ccba{n-2}$, $1\le l\le k-1$.

In fact, if the bottom row  of $D_1$ contains a horizontal arc $\{i,
j\}$, which is different from $\{n-1, n\}$ and if there are $t$ dots
at the left endpoint $i$ of $\{i, j\}$, then we can find another
horizontal arc $\{i', j'\}$ at the top row  of $D_1$ such that there
are $s$ dots at the left endpoint $i'$ of $\{i', j'\}$. Using
vertical arcs $\{i, i'\}$ and $\{j, j'\}$ instead of the horizontal
arcs $\{i, j\}$ and $\{i', j'\}$ in $D_1$, we get another dotted
Brauer diagram $\tilde D_1$. We have $D_1=\tilde D_1 e_{i, j}^{s,
t}$. Note that the number of horizontal arcs in $top (\tilde D_1)$
is $k-1$ if the number of horizontal arcs in $top (D_1)$ is $k$.
Using this method repeatedly, we have $D_1\otimes D_2=D_1'\otimes
e_{i_1,j_1}^{s_1,t_1}\cdots
e_{i_{k-1},j_{k-1}}^{s_{k-1},t_{k-1}}D_2$.

Since  $D_2\in\iocba{n-2}{{k-1}}$, the number of the horizontal arcs
in the top row of the composite of $e_{i_1,j_1}^{s_1,t_1}\cdots
e_{i_{k-1},j_{k-1}}^{s_{k-1},t_{k-1}}$ and $D_2$ is at least $k-1$.
If it is bigger than $k$, then $e_{i_1,j_1}^{s_1,t_1}\cdots
e_{i_{k-1},j_{k-1}}^{s_{k-1},t_{k-1}}D_2=0$ in $\iocba{n-2}{{k-1}}$.
If one loop occurs in the composite of $e_{i_1,j_1}^{s_1,t_1}\cdots
e_{i_{k-1},j_{k-1}}^{s_{k-1},t_{k-1}}$ and $D_2$,
$e_{i_1,j_1}^{s_1,t_1}\cdots
e_{i_{k-1},j_{k-1}}^{s_{k-1},t_{k-1}}D_2=0$ since $\delta_i=0$,
$0\le i\le m-1$. We have $D_1\otimes D_2=0$. In the remaining case,
$e_{i_1,j_1}^{s_1,t_1}\cdots
e_{i_{k-1},j_{k-1}}^{s_{k-1},t_{k-1}}D_2=w\cdot e_{n-3}e_{n-5}\cdots
e_{n-2k+1}$ for some $w\in W_{m, n-2}$. Note that $e_{n-1}w=w
e_{n-1}$,
$$D_1\otimes
D_2=D_1'\otimes e_{i_1,j_1}^{s_1,t_1}\cdots
e_{i_{k-1},j_{k-1}}^{s_{k-1},t_{k-1}}D_2=D_1' w\otimes
e_{n-3}e_{n-5}\cdots e_{n-2k+1}.$$

Since $\{n-1, n\}$ is the unique horizontal arc at the bottom row of
$D_1'$, $D_1' w=w_1 e_{n-1}$ for some $w_1\in W_{m, n}$. Hence,
$D_1\otimes D_2=w_1 e_{n-1}\otimes e_{n-3}e_{n-5}\cdots e_{n-2k+1}$.
We can identify $D_1\otimes D_2$ with $w_1 e_{n-1}e_{n-3}\cdots
e_{n-2k+1}\in \iocba{n}{{k}}$ and vice versa. This proves
$\text{dim}_{F}U_0=\text{dim}_{F}\iocba{n}{{k}}$, where
$U_0=\jcba{n}\otimes_{\ccba{n-2}}\iocba{n-2}{k-1}$.

On the other hand, for any $\alpha\in \iocba{n-2}{k-1}$, let
$\alpha_2^0$ be obtained from $\alpha_2$ by adding  two  vertical
arcs $\{n-1,n-1\}$ and $\{n,n\}$. The $F$-linear map $\phi:
U_0\rightarrow \iocba{n}{k}$ sending $\alpha_1\otimes \alpha_2$ to
$\alpha_1\cdot \alpha_2^0$ is surjective. Since
$\text{dim}_{F}U_0=\text{dim}_{F}\iocba{n}{k}$, it muse be
injective. By the definition of the product of two dotted Brauer
diagrams in \cite{RY}, we can verify that $\phi$ is a $(\ccba{n},
W_{m, n-2k})$-homomorphism. This completes the proof of the claim.

By \cite[(6.3)]{RY}, $\iocba{n}{k} \otimes_{W_{m, n-2k}}
S^{\lambda}\cong S^{(k, \lambda)}$. Therefore,
$$\begin{aligned}
G(S^{(k-1, \lambda)})& = \ccba{n} e_{n-1}\otimes_{\ccba{n-2}}
(\iocba{n-2}{k-1}\otimes_{W_{m, n-2k}} S^{\lambda})\\
&=(\ccba{n} e_{n-1}\otimes_{\ccba{n-2}}
\iocba{n-2}{k-1})\otimes_{W_{m, n-2k}} S^{\lambda}\\
&\cong \iocba{n}{k} \otimes_{W_{m, n-2k}} S^{\lambda}\cong S^{(k,
\lambda)}\\
\end{aligned}$$
\end{proof}

The following theorem is Theorem~\ref{RANK}(b).

\begin{Theorem}
If $n\ge 2$, then $\ccba{n}$ is not (split) semisimple over $F$.
\end{Theorem}

\begin{proof} First, we assume that $n$ is even. A direct
computation shows that the Gram matrix $\G_{\frac n2, 0}$ with
respect to $\Delta(\frac n2, 0)$ is  zero. In particular, $\det
\G_{\frac n2, 0}=0$. By \cite[3.8]{GL}, $\ccba{n}$ is not (split)
semisimple. Suppose $n$ is odd. We have   $n\geq 3$. By
Lemma~\ref{cba3}, there is a non-zero $\phi\in
\text{Hom}_{\ccba{3}}(\Delta(k, \mu), \Delta(1, \lambda))$, where
$\lambda=((1), 0, 0, \cdots, 0)\in \Lambda_m^+(1)$. Write $n=3+2l$
for some $l\in \mathbb N$.
 Applying Proposition~\ref{hom} $l$ times, we get a non-zero homomorphism
from $\Delta(k+l, \mu)$ to  $\Delta(1+l, \lambda)$. By
\cite[3.8]{GL},  $\ccba{n}$ is not (split) semisimple.
\end{proof}


In order to complete the proof of Theorem~\ref{RANK}(a), we need
verify $\tilde \rho_{m, n}=\tilde {\mathbb Z}_{m, n}$.

\begin{Prop} \label{fin} Suppose  $m, n\in \mathbb N$ with $ n\ge
2$.
\begin{itemize} \item[(1)]
$\tilde \rho_{2, n}=\tilde{\rho}_{1, n}=\{k\in \mathbb Z\mid
3-n\le k\le n-3\}\cup \{2k-3\mid 3\le k\le n, k\in \mathbb Z\}$.

\item[(2)] $\tilde{\rho}_{m,n}=\tilde{\rho}_{1, n}\cup \{2-n, n-2\}$ if  $m\geq3$.
\end{itemize}
\end{Prop}

\begin{proof} First, we assume  $m=2$. If $\mu\in \La_m^+(n-2)$ and
$\lambda\in \mathscr A(\mu)$, then  either
$\lambda^{(1)}=\mu^{(1)}$, $\lambda^{(2)}\in \mathscr A(\mu^{(2)})$
or $\lambda^{(2)}=\mu^{(2)}$, $\lambda^{(1)}\in \mathscr
A(\mu^{(1)})$. We can assume $\lambda^{(1)}\in \mathscr
A(\mu^{(1)})$ without loss of generality. Suppose $|\mu^{(1)}|=k$.
Then $k$ can be any integer between $0$ and $n-2$. If $r=\sum_{p\in
Y(\lambda/\mu)} c(p)$, then $r\in \tilde{\rho}_{1, n}$, forcing
$\tilde {\rho}_{2, n}\subset\tilde{\rho }_{1, n}$. Identifying
$\lambda\in \La^+(k)$ with bipartition $(\lambda, (n-k))$, we have
$\tilde {\rho }_{2, n}\supset\tilde{\rho}_{1, n}$. This proves the
first equality in (1). Following \cite{R}, we define
$$ \mathbb Z(n)=\left\{r\in \mathbb Z\mid r= 1-\sum_{p\in
Y(\la/\mu)} c(p), \la\in \La^+(k), \mu\in \La^+(k-2), 2\le k\le
n\right\},
$$ where two boxes in $Y(\lambda/\mu)$ are not in the same column.
By \cite[2.4]{RS},
$$\mathbb Z(n)= \left\{i\in \mathbb Z\mid 4-2n\le i\le n-2\right\}\setminus
\left\{i\in \mathbb Z\mid 4-2n<i\le 3-n, 2\nmid i\right\}.$$
Therefore, the second equality in (1) follows.

Suppose  $m\ge 3$.   If  $\mu$ and $\lambda$ satisfy  one of the
conditions in  Theorem~\ref{admissible}(1) and
Theorem~\ref{admissible}(4), then $\sum_{p\in Y(\la/\mu)} c(p)\in
\tilde{\rho}_{1, n}$. If $\mu$ and $\lambda$ satisfy the conditions
(2) or (3) in Theorem~\ref{admissible}, then there is an $i$, such
that $\mu^{(i)}\rightarrow\lambda^{(i)}$ and
$\mu^{(m-i)}\rightarrow\lambda^{(m-i)}$. In this situation,
$\sum_{p\in Y(\lambda/\mu} c(p)\in \mathfrak T_a+\mathfrak T_b$ with
$|\mu^{(i)}|=a$ and $|\mu^{(m-i)}|=b$, where
\begin{itemize}
\item $\mathfrak T_a=\set {\sum_{p\in Y(\lambda/\mu)} c(p)\mid  \mu\in
\La^+(a), \mu\rightarrow \lambda}$, and
\item $\mathfrak T_a+\mathfrak T_b=\set{i\mid i=x+y, x\in \mathfrak T_a, y\in \mathfrak T_b}$.
\end{itemize}
Note that we can choose a suitable $\mu$ such that $a+b=i$ for all
$i, 0\le i\le n-2$. We claim $$\label{T} \mathfrak T_a=\begin{cases}
\{0\}, &\text{if $a=0$,}\\
 \{i\in \mathbb Z\mid -a\leq i\leq
a\}\setminus \{0\}, &\text{if $a=1, 2$},\\
\{i\in \mathbb Z\mid -a\leq i\leq a\}, & \text{otherwise}\\
\end{cases}
$$
In fact, one can verify the above result directly when  $a\in \{0,1,
2, 3\}$.

Suppose  $\mu\in \La^+(k+1)$ and $\mu\rightarrow\lambda$. If
$\lambda$ has at least two removable nodes, then we can find a box
$q$ which is a removable node  for both $\lambda$ and $\mu$. Let
$\tilde\lambda $(resp. $\tilde\mu$) be obtained from $\lambda$
(resp. $\mu$) by removing $q$. Then  $$\sum_{p\in Y(\lambda/\mu)}
c(p)=\sum_{p\in Y(\tilde\lambda/\tilde\mu)} c(p)\in \mathfrak
T_k=\{-k\leq i\leq k\},$$ the last equality follows from the
induction assumption.

If $\lambda$ has a unique removable node,  then $\lambda=(\lambda_1,
\cdots, \lambda_r)$ with $\lambda_i=\lambda_j$, $1\le i, j\le r$. We
have $\sum_{p\in Y(\lambda/\mu)} c(p)=\lambda_1-r$. Note that
$-1-k\le \lambda_1-r\le k+1$. In any case, we have $\mathfrak
T_{k+1}\subset \{i\in Z\mid-1-k\leq i\leq 1+k\}$.

Conversely,  by the induction assumption, we can write $i=\sum_{p\in
Y(\lambda/\mu)} c(p)$, for some $\lambda\in \La^+(k+1)$ and
$\mu\rightarrow \lambda$ if $-k\le i\le k$. Since any Young diagram
of a partition has at least two addable nodes,  we can choose an
addable node  $q$ for both   $\lambda$ and $\mu$ such that $q$ and
$\lambda/\mu$ are not in the same row. In other words, $i\in
\mathfrak T_{k+1}$. We have

 \begin{itemize}\item $\sum_{p\in Y(\lambda/\mu)} c(p)=-(k+1)$ if
$\lambda=(1, \cdots, 1)\in \La^+(k+2)$ and $\mu=(1, \cdots, 1)\in
\La^+(k+1)$.
\item $\sum_{p\in Y(\lambda/\mu)}
c(p)=k+1$ if $\lambda=(k+2)$ and $\mu=(k+1)$.
\end{itemize}
 Consequently, $\mathfrak T_{k+1}\supset
\{i\in Z\mid-k-1\leq i\leq k+1\}$. This completes the proof of the
claim. Therefore,
$$\cup_{0\le a+b\le n-2} T_a+T_b=\{i\in \mathbb Z\mid 2-n\le i\le
n-2\}.$$ Note that  $i\in \tilde{ \rho}_{1, n}$ if $3-n\le i\le
n-3$. (2) follows immediately.\end{proof}

\textit {Proof of Theorem~\ref{RANK}(a) and (c):}
Theorem~\ref{RANK}(a) follows from Theorem~4.4 and
Proposition~\ref{fin}. Theorem~\ref{RANK}(c) follows from Maschke's
theorem.

\end{document}